\documentclass[12pt]{article}
\usepackage{enumerate}
\usepackage{amssymb,a4wide,latexsym,makeidx,epsfig,fleqn}
\usepackage{amsthm}
\usepackage{amsmath}
\usepackage{enumerate}
\usepackage{float}
\usepackage{color}
\usepackage{amsmath}
\allowdisplaybreaks[4]
\newtheorem{theorem}{Theorem}[section]

\newtheorem{definition}[theorem]{Definition}
\newtheorem{lemma}[theorem]{Lemma}

\begin{document}
\textwidth 150mm \textheight 225mm
\title{The bounds of the spectral radius of general hypergraphs in terms of clique number
\thanks{Supported by the National Natural Science Foundation of China (No. 11871398) and the Seed Foundation of Innovation and Creation for Graduate Students in Northwestern Polytechnical University (No. CX2020190).}}
\author{{Cunxiang Duan$^{a,b}$, Ligong Wang$^{a,b,}$\footnote{Corresponding author.} }\\
{\small $^{a}$School of Mathematics and Statistics, Northwestern
Polytechnical University,} \\ {\small Xi'an, Shaanxi 710129, P.R. China.} \\ {\small $^{b}$Xi'an-Budapest Joint Research Center for Combinatorics, Northwestern Polytechnical University,}
\\ {\small Xi'an, Shaanxi 710129, P.R. China.}
 \\{\small E-mail: cxduanmath@163.com; lgwangmath@163.com}\\}
\date{}
\maketitle
\begin{center}
\begin{minipage}{120mm}
\vskip 0.3cm
\begin{center}
{\small {\bf Abstract}}
\end{center}
{\small The spectral radius (or the signless Laplacian spectral radius) of a general hypergraph is the maximum modulus of the eigenvalues of its adjacency (or  its signless Laplacian) tensor. In this paper, we firstly obtain a lower bound of the spectral radius (or the signless Laplacian spectral radius) of general hypergraphs in terms of clique number. Moreover, we present a relation between a homogeneous polynomial and the clique number of general hypergraphs. As an application, we finally obtain an upper bound of the spectral radius of general hypergraphs in terms of clique number.

}
\end{minipage}
\end{center}

\section{Introduction}
Let $G$ be a general hypergraph with vertex set $V=V(G)=\left\{v_{1}, v_{2}, \ldots, v_{n}\right\}$ and edge set $E=E(G)=\left\{e_{1}, e_{2}, \dots, e_{m}\right\},$ where $e_{1}, e_{2},\ldots, e_{m}\subseteq V(G).$ If $V'\subseteq V$ and $E'\subseteq E$, then $H = (V', E')$ is called a subhypergraph of $G.$ The rank (resp., co-rank) of $G$ is $rank(G)=\max \{|e|\colon\ e \in E\}$ (resp., $corank(G)=\min\{|e|\colon\ e \in E\}$).  If the rank and the co-rank of a hypergraph $G$ are equal to $k$, then $G$ is called $k$-uniform. Obviously, a graph is a $2$-uniform hypergraph. For two vertices $u, v \in V,$ if there exists an edge $e$ such that $u, v \in e \in E,$ then vertices $u$ and $v$ are said to be adjacent. Otherwise, two vertices $u$ and $v$ are said to be nonadjacent. We use $E_{v}$ to denote the set of edges containing $v$ of $G.$ The degree of a vertex $v$ in $G$ is $d_{v} = |E_{v}|.$ In addition, we use $R(v)$ to denote the multiset of edge types in $E_{v}$ for a vertex $v.$

Let $G = (V, E)$  be a general hypergraph. The set $R=\{|e|\colon\ e \in E\}$ is called the set of edge types of $G.$ Then $G$ is also called an $R$-graph. For a set $S$ and a positive integer $i$, let $\binom{S}{i}=\{T\subseteq S \colon\ |T|=i\}.$ A complete $R$-graph on $n$ vertices is an $R$-graph with vertex set $[n]$ and edge set $\cup_{i\in R}\binom{[n]}{i},$ where $[n]=\{1,2,\ldots,n\}$. A clique of a $R$-graph $G$ is a complete $R$-subgraph of $G$. A maximal clique is a clique that cannot extended to a larger clique, and a maximal clique is a clique which covers as many vertices as possible. The clique number $\omega(G)$ of a hypergraph $G$ is the number of vertices in a maximum clique of $G$. If $E(G)=\varnothing$, then we define $\omega(G)=1.$

In 2005, the concept of tensor eigenvalues and the spectra of tensors were independently
introduced by Qi \cite{Q1} and Lim \cite{L}. A tensor $\mathcal{A}=(a_{i_{1}i_{2}\ldots i_{k}})$ with order $k$ and dimension $n$ is a multidimensional array, where $a_{i_{1}i_{2}\ldots i_{k}}\in \mathbb{C}$ and $1\leq i_{1},i_{2},\ldots,i_{k}\leq n.$ For an $n$-dimension complex vector $x=(x_{1}, x_{2}, \ldots, x_{n})^{T}$, $\mathcal{A} x^{k-1}$ is defined as an $n$-dimension complex vector whose $i$-th component is the following
$$(\mathcal{A}x^{k-1})_{i}=\sum\limits^{n}_{i_{2},\ldots,i_{k}=1}a_{ii_{2}\ldots i_{k}}x_{i_{2}}\cdots x_{i_{k}}, \mbox ~{for~any}~ i\in[n].$$

Let $x^{[k-1]}=(x_{1}^{k-1}, x_{2}^{k-1}, \ldots, x_{n}^{k-1})^{T}\in \mathbb{C}^{n}$. Then a number $\lambda\in \mathbb{C}$ is called an eigenvalue of the tensor $\mathcal{A}$ if there exists a nonzero vector $x\in \mathbb{C}^{n}$ such that $$\mathcal{A} x^{k-1}=\lambda x^{[k-1]},$$
and in this case, $x$ is called an eigenvector of $\mathcal{A}$ corresponding to the eigenvalue $\lambda.$ The spectral radius of $\mathcal{A}$ is $\rho(\mathcal{A}) = \max\{| \lambda | \colon\ \lambda \mbox{~is~an~eigenvalue~of~} \mathcal{A}\}.$ In particular, if $x$ is real, then $\lambda$ is called an $H$-eigenvalue. If $x\in\mathbb{R}^{n}_{+}$, then $\lambda$ is called an $H^{+}$-eigenvalue. If $x\in \mathbb{R}^{n}_{++},$ then $\lambda$ is called an $H^{++}$-eigenvalue \cite{Q2}.

\noindent\begin{definition}\label{de:c2}(\cite{BCM})
Let $G=(V,E)$ be a general hypergraph with $n$ vertices, $m$ edges and $rank(G)=k.$ The adjacency tensor $\mathcal{A}(G)$ of $G$ is  $$\mathcal{A}(G)=(a_{i_1,i_2,\ldots,i_k}),~~ 1\leq i_1,i_2,\ldots,i_k\leq n.$$
For all edges $e=\{v_{l_1},v_{l_2},\ldots,v_{l_s}\}\in E$ of cardinality $s\leq k,$
\begin{align*}
a_{i_1,i_2,\ldots,i_k}=\frac{s}{\alpha(s)},  \mbox ~~{where}~ \alpha(s)=\sum_{\substack{k_1+\cdots+k_s=k,\\ k_1,\ldots,k_s\geq 1}}\frac{k!}{k_1!k_2!\cdots k_s!},
\end{align*}
and $i_1,i_2,\ldots,i_k$ are chosen in all possible ways from $L=\{l_1,\ldots,l_s\}$ with at least once for each element of the set $L,$ while each $k_{i}$ represents the times that $l_{i}$ appears in $\{ i_1,i_2,\ldots,i_k\}$. The other positions of the tensor are zero.
\end{definition}

The eigenvalues of a general hypergraph $G$ refer to the eigenvalues of the adjacency tensor of $G.$ And we use $\rho(G)$ to denote the spectral radius of a general hypergraph $G.$

Let $G$ be a general hypergraph with $n$ vertices and $rank(G)=k$, and $x=(x_{1}, x_{2}, \ldots, x_{n})^{T}$ be a vector. For an edge $e=\{i_{1},i_2,\ldots,i_s\}\in E,$ we have

\begin{align*}
x^{T}(\mathcal{A}(e)x^{k-1})=\frac{s}{\alpha(s)}\sum_{\substack{ k_1+\cdots+k_s=k,\\k_1,\ldots,k_s\geq 1}}\frac{k!}{k_1!k_2!\cdots k_s!}x_{i_1}^{k_1}x_{i_2}^{k_2}\cdots x_{i_s}^{k_s}.
\end{align*}
And we also have

\begin{align*}
x^{T}(\mathcal{A}(G)x^{k-1})=\sum_{e\in E(G)}x^{T}(\mathcal{A}(e)x^{k-1}).
\end{align*}

Let $\mathcal{D}(G)$ be a degree diagonal tensor with order $k$ and dimensional $n$, and its diagonal element $d_{ii\ldots i}$ (or simply $d_{i}$) be the degree of a vertex $v_{i}$ of $G$, for all $i \in \{1,2, \ldots, n\}$. Then $\mathcal{Q}(G)=\mathcal{D}(G)+\mathcal{A}(G)$ is the signless Laplacian tensor of the hypergraph $G$. The signless Laplacian eigenvalues of a hypergraph $G$ refer to the eigenvalues of the signless Laplacian tensor of $G$. We use $q(G)$ to denote the signless Laplacian spectral radius of $G$.

Similarly, Let $G$ be a hypergraph with $n$ vertices and $rank(G)=k$, and $x=(x_{1}, x_{2}, \ldots, x_{n})^{T}$ be a vector. For an edge $e=\{i_{1},i_2,\ldots,i_s\}\in E,$ we have

\begin{align*}
x^{T}(\mathcal{Q}(e)x^{k-1})=\sum_{j=1}^{s}x_{i_{j}}^{k}+\frac{s}{\alpha(s)}\sum_{\substack{ k_1+\cdots+k_s=k,\\k_1,\ldots,k_s\geq 1}}\frac{k!}{k_1!k_2!\cdots k_s!}x_{i_1}^{k_1}x_{i_2}^{k_2}\cdots x_{i_s}^{k_s}.
\end{align*}
And we also have

\begin{align*}
x^{T}(\mathcal{Q}(G)x^{k-1})=\sum_{e\in E}x^{T}(\mathcal{Q}(e)x^{k-1}).
\end{align*}

In 1965, Motzkin and Straus \cite{MS} defined a homogeneous polynomial $L(G, x)$ of a graph $G.$

\noindent\begin{definition}\label{de:c3}(\cite{MS})
Let $G=(V,E)$ be a graph with $n$ vertices, and let $x=(x_{1}, x_{2}, \ldots, x_{n})$ be a nonnegative real vector such that $\sum_{i=1}^{n}x_{i}=1$. Then $$L(G,x)=\sum_{\{i,j\}\in E}x_{i}x_{j}.$$
\end{definition}

Moreover, Motzkin and Straus \cite{MS} established a remarkable connection between the clique
number and a homogeneous polynomial $L(G)=\max\limits_{x\in S}L(G,x)$ of a graph $G,$ where $S=\left\{x \colon\ x \in R^{n}_{+}, \sum_{i=1}^{n} x_{i}=1\right\}.$

In 2020, Hou et al. \cite{HCZ} first defined a homogeneous polynomial $L(G, x)$ for a general hypergraph $G$, which generalized the definition $L(G, x)$ from a graph $G$ to a general hypergraph $G.$

\noindent\begin{definition}\label{de:c4}(\cite{HCZ})
Let $G=(V,E)$ be a general hypergraph with $n$ vertices and $rank(G)=k,$ and let $x=(x_{1}, x_{2}, \ldots, x_{n})$ be a nonnegative real vector such that $\sum_{i=1}^{n}x_{i}=1$. Then $$L(G,x)=\sum_{e\in E}\frac{1}{\alpha(s)}\sum_{\substack{ k_1+\cdots+k_s=k,\\k_1,\ldots,k_s\geq 1}}\frac{k!}{k_1!k_2!\cdots k_s!}x_{i_1}^{k_1}x_{i_2}^{k_2}\cdots x_{i_s}^{k_s},$$
where $s$ is the number of vertices of $e$ and $\alpha(s)=\sum\limits_{\substack{k_1+\cdots+k_s=k,\\ k_1,\ldots,k_s\geq 1 }}\frac{k!}{k_1!k_2!\cdots k_s!}.$
\end{definition}

Similarly, for a general hypergraph $G,$  $$L(G)=\max_{x\in S}L(G,x),$$ where $S=\left\{x \colon\ x \in R^{n}_{+}, \sum_{i=1}^{n} x_{i}=1\right\}.$

In 1986, Wilf \cite{W} gave a lower bound on the spectral radius of a graph in terms of clique number, which was inspired by a result of Motzkin and Straus \cite{MS}. In 2007, Bollob\'{a}s and Nikiforov \cite{BN} obtained a number of relations between the number of cliques of a graph $G$ and the spectral radius. In 2007, Lu et al. \cite{LLT} presented some lower and upper bounds for the independence number and the clique number involving the Laplacian eigenvalues of a graph $G.$ In 2008, Liu and Liu \cite{LL} obtained some lower and upper bounds for the independence number and the clique number involving the signless Laplacian eigenvalues of a graph $G.$ In 2009, Rota and Pelillo \cite{RP} gave some new upper and lower bounds on the clique number of graphs. In 2009, Nikiforov \cite{N} gave some new bounds for the clique and independence numbers of a graph in terms of its eigenvalues. In 2013, He et al. \cite{HJZ} gave the bounds on the signless Laplacian spectral radius of graphs in terms of clique number, which disprove the two conjectures on the signless Laplacian spectral radius in \cite{HL}.

Recently, spectral hypergraph theory develops rapidly. There are many work about the
spectral theory of hypergraphs \cite{DW, DWL, DWXL, LZM, XW1, XW2}. However, there are still few studies the relation bewteen the clique number and the (signless Laplacian) spectral radius of hypergraphs. In 2015, Xie and Qi \cite{XQ} mainly researched some inequality relations between the signless (Laplacian) $H$-eigenvalues and the clique (coclique) numbers of uniform hypergraphs. In 2020, Hou et al. \cite{HCZ} first defined a homogeneous polynomial for a general hypergraph, and gave a Motzkin-Straus type result for $\{k-1, k\}$-graphs. And they gave some lower and upper bounds on the spectral radius of $\{k-1, k\}$-graphs in terms of clique number. In this paper, motivated by \cite{HCZ} and \cite{XQ}, we mainly generalize the results of Hou et al. \cite{HCZ} form $\{k-1,k\}$-graphs to $R$-graphs.

This paper is organized as follows. In Section 2, some necessary lemmas and inequalities are given. In Section 3, we give a
bound on the signless Laplacian spectral radius of general hypergraphs in terms of clique number. We also present a relation between a homogeneous polynomial and the clique number of general hypergraphs. Moreover, we also obtain a lower and a upper bound on the spectral radius of general hypergraphs in terms of clique number, respectively.

\section{Preliminaries}
In this section, we mainly give some useful lemmas and inequalities.

\noindent\begin{lemma}\label{le:ch-0}(Maclaurin's inequality \cite{HLP})
Let $x_{1}, x_{2}, \ldots , x_{n}$ be positive real numbers, and $$S_{k}=\frac{\sum\limits_{1 \leq i_{1}<i_{2}<\cdots<i_{k} \leq n} x_{i_{1}} x_{i_{2}} \cdots x_{i_{k}}}{\binom{n}{k}}, \text { for } 1 \leq k \leq n.$$
Then $$S_{1} \geq \sqrt{S_{2}} \geq \sqrt[3]{S_{3}} \geq \cdots \geq \sqrt[n]{S_{n}},$$
and the equality holds if and only if $x_{1} = x_{2} = \cdots = x_{n}.$
\end{lemma}

By Lemma \ref{le:ch-0}, we have $S_{1} \geq \sqrt[k]{S_{k}}$ for any $k \in [n].$ That is, $$\sum_{1\leq i_{1}<i_{2}<\cdots<i_{k} \leq n} x_{i_{1}} x_{i_{2}} \cdots x_{i_{k}} \leq\binom{n}{k}\left(\frac{x_{1}+x_{2}+\cdots+x_{n}}{n}\right)^{k}.$$

\noindent\begin{lemma}\label{le:ch-1}(\cite{ZLKB})
(1). Let $G$ be a general hypergraph, then $\rho(G)$ is an $H^{+}$-eigenvalue of $G.$

(2). If $G$ is connected, then $\rho(G)$ is the unique $H^{++}$-eigenvalue of $G,$ with the unique eigenvector $x \in R^{n}_{++},$ up to a positive scaling coefficient.
\end{lemma}

If $G$ is a connected hypergraph, by Lemma \ref{le:ch-1}, then there exists a unique positive eigenvector $x=(x_{1}, x_{2}, \ldots, x_{n})^{T}$ corresponding to $\rho(G)$ such that $\sum_{i=1}^{n}x_{i}^{k}=1$, where $x$ is called the principal eigenvector of $\mathcal{A}(G)$. For a principal eigenvector $x$ of $\mathcal{A}(G),$ we have $$\rho(G)=x^{T} (\mathcal{A}(G) x^{k-1})=\sum_{e \in E(G)} x^{T} (\mathcal{A}(e) x^{k-1}).$$

Note that $\mathcal{Q}(G)$ is a weakly irreducible tensor. By the Perron-Frobenius theorem for nonnegative tensors \cite{FTP}, there exists a unique positive eigenvector $x=(x_{1}, x_{2}, \ldots, x_{n})^{T}$ corresponding to $q(G)$ such that $\sum_{i=1}^{n}x_{i}^{k}=1$, where $x$ is called the principal eigenvector of $\mathcal{Q}(G)$. Similarly, for a principal eigenvector $x$ of $\mathcal{Q}(G),$ we have $$q(G)=x^{T} (\mathcal{Q}(G) x^{k-1})=\sum_{e \in E(G)} x^{T} (\mathcal{Q}(e) x^{k-1}).$$

\noindent\begin{lemma}\label{le:ch-2}(\cite{HCZ})
Let $G$ be an $R$-graph with clique number $\omega.$ Then
$$\rho(G) \geq \sum_{s \in R}\binom{\omega-1}{s-1},$$
and the equality holds if and only if $G$ is a complete $R$-graph.
\end{lemma}

\section{ Main results}

The following Theorem \ref{th:c1} presents a lower bound on the signless Laplacian spectral radius of a general hypergraph in terms of clique number. This result generalizes Theorem 3.1 of Xie and Qi \cite{XQ} from uniform hypergraphs to general hypergraphs.

\noindent\begin{theorem}\label{th:c1}
Let $G$ be an $R$-graph with clique number $\omega.$ Then $$q(G) \geq 2\sum_{s \in R}\binom{\omega-1}{s-1},$$
and the equality holds if and only if $G$ is a complete $R$-graph.
\end{theorem}

\noindent\textbf{Proof.}
Assume that $G'$ is the maximum complete $R$-subgraph of $G$ with $\omega$ vertices and
$rank(G) = k.$ Let $x$ be a vector such that $x_{i}=\frac{1}{\sqrt[k]{\omega}}$ for $i \in V(G^{\prime})$ and $x_{i}=0$ for otherwise. It is obvious that $\|x\|_{k}^{k}=1.$  By the definition of the signless Laplacian spectral radius of general hypergraphs, we have

\begin{align}\label{eq:1}
q(G) & \geq x^{T} (\mathcal{Q(G)} x^{k-1}) \cr
& =\sum\limits_{\substack{e=\left\{i_{1}, i_{2}, \ldots, i_{s}\right\} \in E\left(G^{\prime}\right)\\s\in R}}(x_{i_{1}}^{k}+x_{i_{2}}^{k}+\cdots+x_{i_{s}}^{k}+\frac{s}{\alpha(s)} \sum\limits_{\substack{k_1+\cdots+k_s=k,\\ k_1,\ldots,k_s\geq 1}} \frac{k !}{k_{1}! k_{2} ! \cdots k_{s}} x_{i_{1}}^{k_{1}} x_{i_{2}}^{k_{2}} \cdots x_{i_{s}}^{k_{s}}) \cr& =\sum\limits_{\substack{e=\left\{i_{1}, i_{2}, \ldots, i_{s}\right\} \in E\left(G^{\prime}\right)\\s\in R}}(\frac{1}{\omega}+\frac{1}{\omega}+\cdots+\frac{1}{\omega}+\frac{1}{\omega} \frac{s}{\alpha(s)} \sum\limits_{\substack{k_1+\cdots+k_s=k,\\ k_1,\ldots,k_s\geq 1}} \frac{k !}{k_{1} ! k_{2} ! \cdots k_{s} !}) \cr
& =\sum\limits_{\substack{e=\left\{i_{1}, i_{2}, \ldots, i_{s}\right\} \in E\left(G^{\prime}\right)\\s\in R}}(\frac{s}{\omega}+\frac{1}{\omega} \frac{s}{\alpha(s)} \alpha(s)) \cr
& =\sum\limits_{s \in R} \sum\limits_{e \in E\left(G^{\prime}\right) \atop |e|=s} \frac{2 s}{\omega}=\sum\limits_{s \in R} \frac{2 s}{\omega}\binom{\omega}{s} \cr
&=2\sum\limits_{s \in R}\binom{\omega-1}{s-1}.
\end{align}

If $q(G)=2 \sum\limits_{s \in R}\binom{\omega-1}{s-1},$  by Inequality (\ref{eq:1}), then we have $q(G)= x^{T} \mathcal{Q}(G) x.$ Since Lemma \ref{le:ch-1} holds for $q(G)$ of a general hypergraph $G$, we know
$x$ is a positive vector. Therefore, we have $\omega = n.$ That is, $G$ is a complete $R$-graph.

If $G$ is a complete $R$-graph, then we have $\omega = n.$ By the definition of the signless Laplacian spectral radius of general hypergraphs, we have $q(G)=2 \sum\limits_{s \in R}\binom{\omega-1}{s-1}.$ \hfill$\square$

If $G$ be a $k$-uniform hypergraph with clique number $\omega$, by Lemma \ref{le:ch-2}, then we have $\rho(G)\geq \binom{\omega-1}{k-1}.$ Therefore, we only consider $R$-graphs such that $|R|>1$ in the following Theorem \ref{th:c2}.

\noindent\begin{theorem}\label{th:c2}
Let $G$ be an $R$-graph with clique number $\omega$, $rank(G)=k$ and $corank(G)=c.$ Then $$\rho(G) \geq \frac{(\omega-c+1)^{c-1}(\omega-k+1)^{k-c}}{(k-1) !}+\frac{(\omega-c+1)^{c-1}}{(c-1) !}.$$
\end{theorem}

\noindent\textbf{Proof.} By Lemma \ref{le:ch-2}, we have
\begin{equation*}\begin{aligned}
\rho(G) & \geq \sum_{s \in R}\binom{\omega-1}{s-1}\\
& \geq\binom{\omega-1}{k-1}+\binom{\omega-1}{c-1} \\
&=\frac{(\omega-1) !}{(k-1) !(\omega-k) !}+\binom{\omega-1}{c-1}\\
&=\frac{(\omega-1) !}{(c-1) !(\omega-c) !} \frac{(c-1) !(\omega-c) !}{(k-1) !(\omega-k) !}+\binom{\omega-1}{c-1} \\
&=\binom{\omega-1}{c-1} \frac{(c-1) !(\omega-c) !}{(k-1) !(\omega-k) !}+\binom{\omega-1}{c-1} \\& \geq \frac{(\omega-c+1)^{c-1}}{(c-1) !} \frac{(c-1) !(\omega-c) !}{(k-1) !(\omega-k) !}+\frac{(\omega-c+1)^{c-1}}{(c-1) !} \\
&=\frac{(\omega-c+1)^{c-1}}{(k-1) !}(\omega-c) \cdots(\omega-k+1)+\frac{(\omega-c+1)^{c-1}}{(c-1) !} \\
& \geq \frac{(\omega-c+1)^{c-1}(\omega-k+1)^{k-c}}{(k-1) !}+\frac{(\omega-c+1)^{c-1}}{(c-1) !},
\end{aligned}\end{equation*}
where $$\binom{\omega-1}{c-1}=\frac{(\omega-1) !}{(c-1) !(\omega-c) !}=\frac{(\omega-1) \cdots(\omega-c+1)}{(c-1) !} \geq \frac{(\omega-c+1)^{c-1}}{(c-1) !}.$$
\hfill$\square$

{\bf Remark.} If $c=k-1>1$, by Theorem \ref{th:c2}, then we have

\begin{equation*}
\begin{aligned}
\rho(G) & \geq \frac{(\omega-k+2)^{k-2}(\omega-k+1)}{(k-1) !}+\frac{(\omega-k+2)^{k-2}}{(k-2) !} \\
&=\frac{(\omega-k+2)^{k-2}}{(k-2) !}\left(\frac{\omega-k+1}{k-1}+1\right)\\
&=\omega \frac{(\omega-k+2)^{k-2}}{(k-1) !}\\
&>\frac{(\omega-k+2)^{k-1}}{(k-1) !}.
\end{aligned}
\end{equation*}

Thus, we know Theorem \ref{th:c2} improves Theorem 3.2 of Hou et al. \cite{HCZ}.

The following Theorem \ref{th:c3} generalizes Theorem 3.3 of Hou et al. \cite{HCZ} from $\{k-1, k\}$-graphs to $R$-graphs.

\noindent\begin{theorem}\label{th:c3}
Let $G=(V,E)$ be an $R$-graph with clique number $\omega$ and $rank(G)=k.$  If $G$ is a complete $R$-graph, or if there exists two nonadjacent vertices $v$ and $v'$ in $G$ such that $R(v)=R(v'),$ then $$L(G)=\sum_{s \in R} \frac{1}{\omega^{k}}\binom{\omega}{s}.$$
\end{theorem}
\noindent\textbf{Proof.}
Assume that $G^{\prime }$ is the maximum complete $R$-subgraph of $G$ with $\omega$ vertices. Let $x$ be an $n$-dimensional vector such that $x_{i}=\frac{1}{\omega}$ for $i \in V(G^{\prime})$ and $x_{i}=0$ for otherwise. It is obvious that $\|x\|_{1}^{1}=1.$ By the definition of $L(G),$ we have

\begin{equation*}
\begin{aligned}
L(G) & \geq L(G, x) \\
&=\sum\limits_{\substack{e=\left\{i_{1}, i_{2}, \ldots, i_{s}\right\} \in E\left(G^{\prime}\right)\\s\in R}} \frac{1}{\alpha(s)} \sum\limits_{\substack{k_1+\cdots+k_s=k,\\ k_1,\ldots,k_s\geq 1}} \frac{k !}{k_{1} ! k_{2} ! \cdots k_{s} !} x_{i_{1}}^{k_{1}} x_{i_{2}}^{k_{2}} \cdots x_{i_{s}}^{k_{s}} \\
&=\sum\limits_{\substack{e=\left\{i_{1}, i_{2}, \ldots, i_{s}\right\} \in E\left(G^{\prime}\right)\\s\in R}} \frac{1}{\omega^{k}} \frac{1}{\alpha(s)}  \sum\limits_{\substack{k_1+\cdots+k_s=k,\\ k_1,\ldots,k_s\geq 1}} \frac{k !}{k_{1} ! k_{2} ! \cdots k_{s} !} \\
&=\sum\limits_{\substack{e=\left\{i_{1}, i_{2}, \ldots, i_{s}\right\} \in E\left(G^{\prime}\right)\\s\in R}} \frac{1}{\omega^{k}} \frac{1}{\alpha(s)} \alpha(s) \\
&=\sum_{s \in R} \sum_{e \in E\left(G^{\prime}\right) \atop |e|=s} \frac{1}{\omega^{k}} \\
&=\sum_{s \in R} \frac{1}{\omega^{k}}\binom{\omega}{s}.
\end{aligned}
\end{equation*}

In the following, we prove

\begin{align}\label{eq:2}
L(G) \leq \sum_{s \in R} \frac{1}{\omega^{k}}\binom{\omega}{s}
\end{align}
by induction on $n$.

If $n \leq c-1,$ then the edge set of $G$ is empty. It is obvious that $L(G) \leq \sum_{s \in R} \frac{1}{\omega^{k}}\binom{\omega}{s}$ holds. Assume that Inequality (\ref{eq:2}) holds for less than $n$ vertices. Without loss of generality, supposed that $n$ is sufficient large and $x$ is a vector such that $L(G)=L(G, x) .$ And $x_{1} \geq x_{2} \geq \cdots \geq x_{t}>x_{t+1}>\cdots>x_{n}=0.$

If $t<n,$ then we have $x_{n}=0 .$ We can obtain $G_{1}$ from $G$ by deleting the vertex $n$ and the edges containing the vertex $n$. Let $G_{1}$ be a general hypergraph with clique number $\omega_{1}$ and $\omega_{1} \leq \omega .$ By the induction hypothesis, we have

$$L(G)=L\left(G_{1}\right) \leq \sum_{s \in R}\frac{1}{\omega_{1}^{k}}\binom{\omega_{1}}{s} \leq \sum_{s \in R} \frac{1}{\omega^{k}}\binom{\omega}{s},$$
where $\left(\frac{1}{\omega_{1}}\right)^{k}\binom{\omega_{1}}{s}$ is an monotonically increasing function of $\omega_{1}.$

If $t=n,$ then we have $x_{1} \geq x_{2} \geq \cdots \geq x_{n}>0$.

{\bf Case 1.} If $G$ is a complete $R$-graph, then we have $\omega=n .$ By the definition of $L(G)$ and triangle inequality, we have

\begin{align}\label{eq:3}
L(G) &=\max _{x \in S} L(G, x) \cr
&=\max _{x \in S} \sum_{e=\left\{i_{1}, i_{2}, \ldots, i_{s}\right\} \in E(G)\atop s\in R} \frac{1}{\alpha(s)} \sum_{ k_{1}+k_{2}+\cdots+k_{s}=k,\atop k_{1}, \ldots, k_{s} \geq 1} \frac{k !}{k_{1} ! k_{2} ! \cdots k_{s} !} x_{i_{1}}^{k_{1}} x_{i_{2}}^{k_{2}} \cdots x_{i_{s}}^{k_{s}} \cr
& \leq \max _{x \in S} \sum_{e \in E(G)\atop |e|=c} \frac{1}{\alpha(c)} \sum_{k_{1}+k_{2}+\cdots+k_{c}=k,\atop k_{1}, \ldots, k_{c} \geq 1 } \frac{k !}{k_{1} ! k_{2} ! \cdots k_{c} !} x_{i_{1}}^{k_{1}} x_{i_{2}}^{k_{2}} \cdots x_{i_{c}}^{k_{c}} \cr
&+\cdots+\max _{x \in S} \sum_{e \in E(G) \atop|e|=k} x_{i_{1}} x_{i_{2}} \cdots x_{i_{k}}.
\end{align}

Without loss of generality, assume that there exist two integers $d \geq 1$ and $u \geq 0$ such taht $k=s d+u$ for any $s \in R .$ We know that

\begin{align}\label{eq:4}
\alpha(s) &=\sum_{k_{1}, \ldots, k_{s} \geq 1,\atop k_{1}+\cdots+k_{s}=k } \frac{k !}{k_{1} ! k_{2} ! \cdots k_{s} !} \cr
&=a_{1}\binom{s}{1}+a_{2}\binom{s}{2}+a_{3}\binom{s}{2}+\cdots+a_{b}\binom{s}{u} \cr
&=a_{1} \frac{s !}{(s-1) ! 1 !}+a_{2} \frac{s !}{(s-2) ! 2 !}+a_{3} \frac{s !}{(s-2) ! 2 !}+\cdots+a_{b} \frac{s !}{(s-u) ! u !} \cr
&=s !\left( \frac{a_{1}}{(s-1) ! 1 !}+ \frac{a_{2}}{(s-2) !2 !}+ \frac{a_{3}}{(s-2) !2 !}+\cdots+ \frac{a_{b}}{(s-u) ! u !}\right),
\end{align}
where $$a_{1}=\frac{k !}{(k-s+1) ! 1 ! \cdots 1 !},~~a_{2}=\left\{
\begin{array}{ll}
\frac{k !}{(k-s) ! 2 ! 1 ! \cdots 1 !},& \mbox {if}   ~k-s = 2,
\\
2! \frac{k !}{(k-s) ! 2 ! 1 ! \cdots 1 !},& \mbox {if}   ~k-s \neq 2,
\end{array}
\right. ,$$
and $$a_{3}=\left\{
\begin{array}{ll}
\frac{k !}{(k-s-1) ! 3 ! 1 ! \cdots 1 !},& \mbox {if}   ~k-s-1 = 3,
\\
2! \frac{k !}{(k-s-1) ! 3 ! 1 ! \cdots 1 !},& \mbox {if}   ~k-s-1 \neq 3,
\end{array}
\right. ,~~a_{b}=\frac{k !}{\underbrace{(d+1) ! \cdots(d+1) !}_{u~times} d ! \cdots d !}.$$

Since $G$ is a complete $R$-graph, we also have

\begin{align*}
&\sum_{e=\left\{i_{1}, i_{2}, \ldots, i_{s}\right\} \in E(G),\atop |e|=s}\sum_{k_{1}+k_{2}+\cdots+k_{s}=k,\atop k_{1}, \ldots, k_{s} \geq 1} \frac{k !}{k_{1} ! k_{2} ! \cdots k_{s} !} x_{i_{1}}^{k_{1}} x_{i_{2}}^{k_{2}} \cdots x_{i_{s}}^{k_{s}}\\
&=\frac{k !}{(k-s+1) ! 1 ! \cdots 1 !} \sum_{e=\left\{i_{1}, i_{2}, \ldots, i_{s}\right\} \in E(G)} x_{i_{1}} x_{i_{2}} \cdots x_{i_{s}}\left(x_{i_{1}}^{k-s}+\cdots+x_{i_{s}}^{k-s}\right)\\
&+\frac{k !}{(k-s) ! 2 ! 1 ! \cdots 1 !} \sum_{e=\left\{i_{1}, i_{2}, \ldots, i_{s}\right\} \in E(G)} x_{i_{1}} x_{i_{2}} \cdots x_{i_{s}}\left(x_{i_{1}}^{k-s-1} x_{i_{2}}+\cdots+x_{i_{s-1}} x_{i_{s}}^{k-s-1}\right)+\cdots\\
&+\frac{k !}{(d+1) ! \cdots(d+1) ! d ! \cdots d !} \sum_{e=\left\{i_{1}, i_{2}, \dots, i_{s}\right\} \in E(G)} x_{i_{1}}^{d} x_{i_{2}}^{d} \cdots x_{i_{s}}^{d}\left(x_{i_{1}} \cdots x_{i_{u}}+\cdots+x_{i_{s-u+1}} \cdots x_{i_{s}}\right)\\
&=\frac{k !}{(k-s+1) ! }[x_{1}^{k-s+1}\sum_{e=\left\{1, i_{2}, \ldots, i_{s}\right\} \in E(G),\atop 2\leq i_{2}< \cdots < i_{s}\leq n} x_{i_{2}}\cdots x_{i_{s}} +\cdots+x_{n}^{k-s+1}\sum_{e=\left\{i_{1}, i_{2}, \ldots, n\right\} \in E(G),\atop 1\leq i_{1}< \cdots < i_{s-1}\leq n-1} x_{i_{1}}\cdots x_{i_{s-1}} ]\\
&+\frac{k !}{(k-s) ! 2 !}[x_{1}^{k-s} x_{2}^{2}\sum_{e=\left\{1, 2, \ldots, i_{s}\right\} \in E(G),\atop 3\leq i_{3}< \cdots < i_{s}\leq n} x_{i_{3}}\cdots x_{i_{s}} +\cdots+x_{n-1}^{2} x_{n}^{k-s}\sum_{e=\left\{i_{1}, i_{2}, \ldots, n-1,n \right\} \in E(G),\atop 1\leq i_{1}< \cdots < i_{s-2}\leq n-2} x_{i_{2}}\cdots x_{i_{s}} ]\\
&+\cdots+\frac{k !}{(d+1) ! \cdots(d+1) ! d ! \cdots d !}[x_{1}^{d} \cdots x_{s}^{d}\sum_{e=\left\{i_{1}, i_{2}, \ldots,i_{u}, n-s+u+1,\ldots, n\right\} \in E(G), \atop 1\leq i_{1}< \cdots < i_{u}\leq n-s+u} x_{1}\cdots x_{u} \\ &+\cdots+ x_{1}^{d} \cdots x_{s}^{d}\sum_{e=\left\{1, 2, \ldots,s-u, i_{s-u+1},\ldots, i_{s}\right\} \in E(G),\atop s-u+1\leq i_{s-u+1}< \cdots < i_{s}\leq n} x_{i_{s-u+1}}\cdots x_{i_{s}}].
\end{align*}

Thus, by Maclaurin's inequality, we have
\begin{align*}
&\sum_{e=\left\{i_{1}, i_{2}, \ldots, i_{s}\right\} \in E(G),\atop |e|=s}\sum_{k_{1}+k_{2}+\cdots+k_{s}=k,\atop k_{1}, \ldots, k_{s} \geq 1} \frac{k !}{k_{1} ! k_{2} ! \cdots k_{s} !} x_{i_{1}}^{k_{1}} x_{i_{2}}^{k_{2}} \cdots x_{i_{s}}^{k_{s}}\\
&\leq \frac{k !}{(k-s+1) ! } [x_{1}^{k-s+1}\left(\frac{\sum_{2 \leq i \leq n} x_{i}}{n-1}\right)^{s-1}\binom{n-1}{s-1}+\cdots+x_{n}^{k-s+1}\left(\frac{\sum_{1 \leq i \leq n-1} x_{i}}{n-1}\right)^{s-1}\binom{n-1}{s-1}]\\
&+\frac{k !}{(k-s) ! 2 ! }[x_{1}^{k-s} x_{2}^{2}\left(\frac{\sum_{3 \leq i \leq n} x_{i}}{n-2}\right)^{s-2}\binom{n-2}{s-2}+\cdots+x_{n-1}^{2} x_{n}^{k-s}\left(\frac{\sum_{1 \leq i \leq n-2} x_{i}}{n-2}\right)^{s-2}\binom{n-2}{s-2}]\\
&+\cdots+\frac{k !}{(d+1) ! \cdots(d+1) ! d ! \cdots d !}[x_{1}^{d} \cdots x_{s}^{d}\left(\frac{\sum_{1 \leq i \leq n-s+u} x_{i}}{n-s+u}\right)^{u}\binom{n-s+u}{u}\\
&+\cdots+x_{1}^{d} \cdots x_{s}^{d}\left(\frac{\sum_{s-u+1 \leq i \leq n} x_{i}}{n-s+u}\right)^{u}\binom{n-s+u}{u}].
\end{align*}
The equality holds if and only if $x_{1}=x_{2}=\cdots=x_{n}=\frac{1}{n} .$

For any $1\leq i \leq s$, we have $$\binom{n-i}{s-i}\binom{n}{i}=\binom{n}{s}\binom{s}{i}=\frac{n !}{(n-s) ! } \frac{1}{(s-i) !i !}.$$
Let $x_{1}=x_{2}=\cdots=x_{n}=\frac{1}{n}.$ Therefore, we have

\begin{align}\label{eq:5}
&\sum_{e=\left\{i_{1}, i_{2}, \ldots, i_{s}\right\} \in E(G),\atop |e|=s}\sum_{k_{1}+k_{2}+\cdots+k_{s}=k,\atop k_{1}, \ldots, k_{s} \geq 1} \frac{k !}{k_{1} ! k_{2} ! \cdots k_{s} !} x_{i_{1}}^{k_{1}} x_{i_{2}}^{k_{2}} \cdots x_{i_{s}}^{k_{s}} \cr &
\leq a_{1}\left(\frac{1}{n}\right)^{k}\binom{n-1}{s-1}\binom{n}{1}+a_{2}\left(\frac{1}{n}\right)^{k}\binom{n-2}{s-2}\binom{n}{2}+\cdots+a_{b}\left(\frac{1}{n}\right)^{k}\binom{n-s+u}{u}\binom{n}{s-u} \cr &
=\left(\frac{1}{n}\right)^{k} \frac{n !}{(n-s) !}\left( \frac{a_{1}}{(s-1) ! 1 !}+ \frac{a_{2}}{(s-2) !2 !}+\cdots+ \frac{a_{b}}{u !(s-u) !}\right).
\end{align}

By (\ref{eq:4}) and (\ref{eq:5}), we have

\begin{align}\label{eq:6}
&\frac{1}{\alpha(s)} \sum_{e=\left\{i_{1}, i_{2}, \ldots, i_{s}\right\} \in E(G), \atop |e|=s} \sum_{k_{1}+k_{2}+\cdots+k_{s}=k,\atop k_{1}, \ldots, k_{s} \geq 1 }\frac{k !}{k_{1} ! k_{2} ! \cdots k_{s} !} x_{i_{1}}^{k_{1}} x_{i_{2}}^{k_{2}} \cdots x_{i_{s}}^{k_{s}}\cr &\leq\left(\frac{1}{n}\right)^{k} \frac{n !}{(n-s) ! s !}=\left(\frac{1}{n}\right)^{k}\binom{n}{s}.
\end{align}

By (\ref{eq:3}) and (\ref{eq:6}), we have

$$L(G) \leq \sum_{s \in R}\left(\frac{1}{n}\right)^{k}\binom{n}{s}.$$

{\bf Case 2.} If there exist two nonadjacent vertices $v$ and $v'$ in $G$ such that $R(v)=R(v'),$ then we prove that there exists a vector $y$ such that $L(G, y) \geq L(G, x) .$ By the definition of $L(G),$ we denote $L^{v}(G)$ as

$$L^{v}(G)=\sum_{e \in E_{v},|e|=s\in R} \frac{1}{\alpha(s)} \frac{k !}{k_{1} ! k_{2} ! \cdots k_{s} !} x_{i_{1}}^{k_{1}} x_{i_{2}}^{k_{2}} \cdots x_{i_{s}}^{k_{s}},$$
for any a vertex $v \in V(G).$

Let $y$ be a vector such that $y_{i}=x_{i}$ for $i \neq v, v',$ $y_{v}=x_{v}+x_{v'}$ and $y_{v'}=0 .$ And let $z$ be a vector such that $z_{i}=x_{i}$ for $i \neq v$ and $z_{v}=x_{v'} .$ It is obvious that $y, z \in R_{+}^{n} .$ For any positive integer $k,$ we have $\left(x_{v}+x_{v'}\right)^{k} \geq x_{v}^{k}+x_{v'}^{k} .$ Hence, we have $$L^{v}(G,y)\geq L^{v}(G,x)+L^{v}(G,z).$$ Moreover, without loss of generality, assume that $L^{v}(G,z)\geq  L^{v'}(G,x).$ Then we have $$L(G,y)- L(G,x)=L^{v}(G,y)- L^{v}(G,x)-L^{v'}(G,x)\geq L^{v}(G,z)-L^{v'}(G,x)\geq 0.$$

Thus, we can obtain $G_{1}$ from $G$ by deleting the vertex $v'$ and the edges containing $v'.$ Similarly, by the induction hypothesis, we have
$$L(G)=L(G,x)\leq L(G,y)\leq L(G_{1})\leq\sum_{s\in R}(\frac{1}{\omega})^{k}\binom{\omega}{s}.$$

Therefore, we have $$L(G)=\sum_{s \in R} \frac{1}{\omega^{k}}\binom{\omega}{s},$$
where $G$ is a complete $R$-graph, or there exists two nonadjacent vertices $v$ and $v'$ in $G$ such that $R(v) = R(v').$
\hfill$\square$

The following Theorem \ref{th:c4} generalizes Theorem 3.4 of Hou et al. \cite{HCZ} from $\{k-1, k\}$-graphs to $R$-graphs.

\noindent\begin{theorem}\label{th:c4}
Let $G = (V, E)$ be an $R$-graph with clique number $\omega$ and $rank(G) = k.$ If
$G$ is a complete $R$-graph, or if there exists two nonadjacent vertices $v$ and $v'$ in $G$ such that $R(v) = R(v'),$ then $$\rho(G)\leq\sum_{s\in R}k(\frac{U}{\omega})^{k}\binom{\omega}{s},$$
where $U$ is the sum of entries of the principal eigenvector.
\end{theorem}

\noindent\textbf{Proof.}
Let $x$ be the principal eigenvector of $G$ with $n$ vertices, and $y_{i}=\frac{x_{i}}{U},$ for $1\leq i \leq n.$ It is obvious that $\sum_{i=1}^{n}y_{i}=1.$  By Theorem \ref{th:c3} and the definition of $\rho(G),$ we have

\begin{align*}
\frac{\rho(G)}{U^{k}}&=\sum_{e=\left\{i_{1}, i_{2}, \ldots, i_{s}\right\} \in E(G),\atop s\in R}\frac{s}{\alpha(s)}\sum\limits_{\substack{k_1+\cdots+k_s=k,\\ k_1,\ldots,k_s\geq 1}} \frac{k !}{k_{1} ! k_{2} ! \cdots k_{s} !} y_{i_{1}}^{k_{1}} y_{i_{2}}^{k_{2}} \cdots y_{i_{s}}^{k_{s}}\cr
&\leq\sum_{e=\left\{i_{1}, i_{2}, \ldots, i_{s}\right\} \in E(G),\atop s\in R}\frac{k}{\alpha(s)}\sum\limits_{\substack{k_1+\cdots+k_s=k,\\ k_1,\ldots,k_s\geq 1}} \frac{k !}{k_{1} ! k_{2} ! \cdots k_{s} !} y_{i_{1}}^{k_{1}} y_{i_{2}}^{k_{2}} \cdots y_{i_{s}}^{k_{s}}\cr
&\leq \sum_{s\in R}k\frac{1}{\omega^{k}}\binom{\omega}{s}.
\end{align*}

Therefore, we have
$$\rho(G)\leq\sum_{s\in R}k(\frac{U}{\omega})^{k}\binom{\omega}{s}.$$
\hfill$\square$

\end{document}